\documentclass[a4paper,12pt]{article}
\usepackage[estonian]{babel}
\usepackage[T1]{fontenc}
\usepackage[latin1]{inputenc}
\usepackage{amsmath}
\usepackage{graphicx}

\usepackage{amsthm}
\usepackage{lineno}

\begin{document}
\linenumbers
\nolinenumbers
\begin{center}

\bf \large
    {A version of the simplex method for solving linear systems of inequalities and linear programming problems}
\end{center}
\begin{center}
\bf\large
    {Jaan Ubi}
\end{center}
\begin{center}
\medskip
\bf \large
               {\bf jaanbi.jb@gmail.com }
\end{center}

 \begin{center}
\bf\large
    {Evald Ubi}
\end{center}

\medskip
{\bf Tallinn Universityy of Tecnology; Evald.Ubi@ttu.ee }

------------------------------------------------------------------------------------------------

 \noindent {\normalsize  In order to find a non-negative solution to a system of inequalities, the corresponding dual problem is composed, which has a suitable unity basic matrix. In such a formulation, the objective function is replaced by set of constraints based on that function. Additional constraints can be used for accelerating calculations in the first phase of the simplex method.  As an example the solution of the Klee-Minty hypercube problem is described in detail.

\vskip 1.5cm
\section{\bf Introduction}

    The elimination method for solving systems of linear equations has been discovered over 2000 years ago. That method remained unknown in Europe until C. Gauss rediscovered it in the 19th century. The simplex algorithm is emerged in the late 1940s from the US mathematician George Dantzig, who had spent the second world war investigating ways to increase the logistical efficiency of the US air force.  The simplex method exhibits exponential growth in the worst case, but its performance in practice has been outstanding.

       Although the theory of linear equations was worked out in details in the 19th century, there was only limited progress in what today is considered a related subject - systems of linear inequalities. The first to deal with these ideas in some detail was J. Fourier (1768-1830) in the 1820s. He was interested in various types of problems in which inequalities appeared (mechanics, probability, statistics). Fourier worked out both algebraic and geometric methods of finding the region of solutions, this was augmented by T.Motzkin. The Fourier-Motzkin elimination method is very elegant theoretically, however the elimination of each variable adds new inequalities to the remaining system, but the number of these new inequalities grows exponentially, see /1/.

       Linear inequalities are more complicated than linear equations. Equations as constraints can be transformed into linear inequalities by replacing each equation with the opposite pair of inequalities, but linear inequalities cannot be transformed into equations. We would first learn to solve the inequalities and then worry abut minimizing over the set of feasible solutions.

    In this work, we consider the calculation of the non-negative solution to the systems of linear inequalities.  We replace the linear programming problem with a certain number of systems of linear inequalities. We compose the dual problem for a system of linear inequalities and solve it using the simplex method.

    Systems of linear inequalities and linear programming problems can also be solved by the method of least squares, see /2,3/. Orthogonal transformations can be used to find the solution of minimum norm to a system of linear inequalities.

\vskip 1.5cm
\section{\bf Introduction}

  Let us consider a system
    \begin{equation}\
\ {Ax \geq b.}
\ x\geq 0,
\end{equation}

where $A$ is a given matrix with dimensions $m\times n$, $b$ is a given vector with dimensions $m\times 1$ and $x$ $n\times 1$.

    Let us take as goal vector 0, which all components are zeroes, 0=(0,...,0) with dimension $1\times n$

 \begin{equation}
  z = (0,x)\rightarrow min,
 \ \,\,\,\,\,{Ax \geq  b,}
 \end{equation}
$$ x\geq 0, $$

and compose the dual problem,

  \begin{equation}\
  w = (y, b)\rightarrow max,
\ {yA +Is = 0,}
\ y, s\geq 0,
\end{equation}
where $I$ is a unity matrix, dual variables $y$ is a $1\times m$ and $s$ $n\times 1$ vectors.

    If problem (1) has no solution, then dual problem (3) is unbounded. Otherwise, $ z_{min} = w_{max} =0$, the w-row gives a solution to the system of inequalities (1).

\vskip 1.5cm
\section{\bf Two numerical examples}

{\bf Example 3.1.}
     $$ x_{1}+3x_{2}\,\,\,\,\geq 2 $$
     $$\,\,\,\,x_{2}\,\,\,\,\geq 1$$
    $$-x_{1}\,\,\,\,+x_{3}\geq -2 $$
              $$ x\geq 0. $$
    We solve the dual problem (3) using the simplex method.

\begin{tabular}{|c|c|c|c|c|c|}
  \hline
  \hline
  $1$ &2 & 3 & 4 & 5&6  \\
  \hline
   2&1&-2&0&0&0 \\
  1&0&-1&1&0&0 \\
  3&1&0&0&1&0 \\
   0&0&1&0&0&1 \\
  \hline

\end{tabular}

    Table 1

\begin{tabular}{|c|c|c|c|c|c|}
  \hline
  \hline
  $1$ &2 & 3 & 4 & 5&6  \\
  \hline
   0&1&0&-2&0&0 \\
   1&0&-1&1&0&0 \\
   0&1&3&-3&1&0 \\
   0&0&1&0&0& 1 \\
  \hline

\end{tabular}

Table 2

\begin{tabular}{|c|c|c|c|c|c|}
  \hline
  \hline
  $1$ &2 & 3 & 4 & 5&6  \\
  \hline
   0&0&-3&1&-1&0 \\
  1&0&-1&1&0&0 \\
  0&1&3&-3&1&0 \\
  0&0&1&0&0&1 \\
  \hline

\end{tabular}

Table 3

\begin{tabular}{|c|c|c|c|c|c|}
  \hline
  \hline
  $1$ &2 & 3 & 4 & 5&6  \\
  \hline
   -1&0&-2&0&-1&0 \\
  1&0&-1&1&0&0 \\
  3&1&0&0&1&0 \\
  0&0&1&0&0&1 \\
  \hline

\end{tabular}

Table 4

    In these tables, the variable leaving the base is specified by the lexicographically minimal row, see /4/.   In w-row are the value of active variable $x_{1}, x_{1}= 2,$ see table 2. In table 3 $x_{1}=-1, x_{2}=1$ and in table 4 the solution to the system of inequalities,$x=(0, 1, 0)^{T},$ see /5/.

    {\bf Remark 3.1.}
    If we multiply the second inequality by three, problem 3.2 is solved in one step.

{\bf Example 3.2.}
    In paper /6/ a cycling example is given.
     $$ z = -2x_{1}-3x_{2}+x_{3} +12x_{4} \rightarrow min $$
     $$ 2x_{1}+9x_{2}-x_{3} -9x_{4}\geq 0 $$
      $$ -x_{1}/3-x_{2}+x_{3}/3 +2x_{4} \geq 0$$
    $$ -2x_{1}-3x_{2}+x_{3} +12x_{4}\geq 18 $$
              $$ x\geq 0. $$
\begin{tabular}{|c|c|c|c|c|c|c|}
  \hline
  \hline
  $1$ &2 & 3 & 4 &5&6&7 \\
  \hline
  0&0&18&0&0&0&0 \\
   2&-0,33&-2&1&0&0&0 \\
  9&-1&-3&0&1&0&0 \\
  -1&0,33&1&0&0&1&0 \\
  -9 &2&12 &0&0&0&1 \\
  \hline
\end{tabular}

  Table 5

\begin{tabular}{|c|c|c|c|c|c|c|}
  \hline
  \hline
  $1$ &2 & 3 & 4 &5&6&7 \\
  \hline
 18&-6&0&0&0&-18&0 \\
   0&0,33&0&1&0&2&0 \\
  6&0&0&0&1&3&0 \\
 -1&0,33&1&0&0&1&0 \\
   3 &-2&0 &0&0&-12&1 \\
  \hline
\end{tabular}

 Table 6

\begin{tabular}{|c|c|c|c|c|c|c|}
  \hline
  \hline
  $1$ &2 & 3 & 4 &5&6&7 \\
  \hline
 0&6&0&0&0&54&-6 \\
   0&0,33&0&1&0&2&0 \\
  0&4&0&0&1&27&-2 \\
 0&-0,33&1&0&0&-3&0.33 \\
   1 &-0,67&0 &0&0&-4&0,33 \\
  \hline
\end{tabular}

  Table 7

\begin{tabular}{|c|c|c|c|c|c|c|}
  \hline
  \hline
  $1$ &2 & 3 & 4 &5&6&7 \\
  \hline
 0&-2&0&0&-2&0&-2 \\
   0&0,037&0&1&-0,074&0&0,148 \\
  0&0,148&0&0&0,037&1&-0,074 \\
 0&0,11&1&0&0,111&0&0,111 \\
   1 &-0,074&0 &0&0,148&0&0,037 \\
  \hline
\end{tabular}

Table 8

The found solution $x*=(0,2,0,2)^{T}$ is one of the optimal solutions to this problem.

\vskip 1.5cm
\section{\bf Systems of inequalities instead of linear programming}

    Studying the system of inequalities could be the first step in composing a linear programming problem.  You can then specify the required resource quantities, specify one or more target functions, or replace one of these functions with a constraint. The objective function of a problem with a large number of variables and constraints can be unbounded if there are no significant substantive constraints. In the second option, the constraints may be contradictory. The impossibility would be more disguised (unless it arose through a simple keying error). Assuming that we are modelling a situation which we know to be realizable, it should be possible to construct a feasible (but probably non-optimal) solution to the problem. Any constraints violated for that solution must have been modelled wrongly. They are too restrictive and should be reconsidered.
 Most package programs will print out the infeasible solution obtained at the point when the program gives up. In solving this problem, we make sure that the presented method is suitable for the addition of constraints.

                                    We describe this process on the Klee-Minty example , which we solve with different constraints for the objective function.

    The replacement of the linear programming problem with a certain number of systems of linear inequalities has been considered in the work /7/.

   {\bf Example 4.1.}

    Let us study the Klee-Minty problem for four variables. The number of steps using simplex method grows exponentially, see /1/.

 $$ z = -8x_{1}-4x_{2}-2x_{3} -x_{4} \rightarrow min $$
     $$ -x_{1} \,\,\,\,\,\,\,\,\,\,\,\,\,\,  \,\,\,\,\,\,\,\,\,\,  \,\,\,\,\,\,\,                   \geq -5 $$
      $$ -4x_{1}-x_{2}\,\,\,\,\,\,\,\,\, \,\,\,\,\,\,\,\,\,\,\,\,\,\,\,\, \geq -25 $$
    $$ -8x_{1}-4x_{2}-x_{3},\,\,\, \,\,\,\,\,\,\,\,\, \geq -125 $$
    $$ -16x_{1}-8x_{2}-4x_{3} -x_{4} \geq -625 $$
              $$ x\geq 0. $$
    The solution is $x=(0, 0, 0, 625)^{T}.$

   Let us solve three system at the same time: $-z \geq 500$,  $-z \geq 600$  and  $-z \geq 700$.

\begin{tabular}{|c|c|c|c|c|c|c|c|c|c|c|}
  \hline
  \hline
  1&2 & 3 & 4 &5&6&7&8& 9&10&11 \\
  \hline
  700&600&500&-5&-25&-125&-625&0&0&0&0 \\
   8&8&8&-1&-4&-8&-16&1&0&0&0 \\
  4&4&4&0&-1&-4&-8&0&1&0&0 \\
  2&2&2&0&0&-1&-4&0&0&1&0 \\
   1&1&1&0&0&0&-1&0&0&0&1 \\
  \hline
  \end{tabular}

  Table 1

\begin{tabular}{|c|c|c|c|c|c|c|c|c|c|c|}
  \hline
  \hline
  1&2 & 3 & 4 &5&6&7&8& 9&10&11 \\
  \hline
  200&100&0&57,5&225&375&375&-62,5&0&0&0 \\
   1&1&1&-0,125&-0,5&-1&-2&0,125&0&0&0 \\
  0&0&0&0,5&1&0&0&-0,5&1&0&0 \\
  0&0&0&0,25&1&1&0&-0,25&0&1&0 \\
  0&0&0&0,125&0,5&1&1&-0,125&0&0&1 \\
  \hline
  \end{tabular}

  Table 2

  \begin{tabular}{|c|c|c|c|c|c|c|c|c|c|c|}
  \hline
  \hline
  1&2 & 3 & 4 &5&6&7&8& 9&10&11 \\
  \hline
  200&100&0&10,625&37,5&0&0&-15,625&0&0&-375 \\
   1&1&1&0,125&0,5&1&0&-0,125&0&0&2 \\
  0&0&0&0,5&1&0&0&-0,5&1&0&0 \\
  0&0&0&0,25&1&1&0&-0,25&0&1&0 \\
  0&0&0&0,125&0,5&1&1&-0,125&0&0&1 \\
  \hline
  \end{tabular}

  Table 3

  \begin{tabular}{|c|c|c|c|c|c|c|c|c|c|c|}
  \hline
  \hline
  1&2 & 3 & 4 &5&6&7&8& 9&10&11 \\
  \hline
  200&100&0&1,25&0&-37,5&0&-6,25&0&-37,5&-375 \\
   1&1&1&0&0&0,5&0,5&0&0&-0,5&2 \\
  0&0&0&0,25&0&-1&0&-0,25&1&-1&0 \\
  0&0&0&0,25&1&1&0&-0,25&0&1&0 \\
  0&0&0&0&0&0,5&1&0&0&-0,5&1 \\
  \hline
  \end{tabular}

  Table 4

  \begin{tabular}{|c|c|c|c|c|c|c|c|c|c|c|}
  \hline
  \hline
  1&2 & 3 & 4 &5&6&7&8& 9&10&11 \\
  \hline
  200&100&0&0&0&-32,5&0&-5&-5&-32,5&-375 \\
   1&1&1&0&0&0,5&0&0&0&-0,5&2 \\
  0&0&0&1&0&-4&0&-1&4&-4&0 \\
  0&0&0&0&1&2&0&0&-1&2&0 \\
  0&0&0&0&0&0,5&1&0&0&-0,5&1 \\
  \hline
  \end{tabular}

  Table 5

 \begin{tabular}{|c|c|c|c|c|c|c|c|c|c|c|}
  \hline
  \hline
  1&2 & 3 & 4 &5&6&7&8& 9&10&11 \\
  \hline
  100&0&-100&0&0&-82,5&0&-5&-5&17,5&-575 \\
   1&1&1&0&0&0,5&0&0&0&-0,5&2 \\
  0&0&0&1&0&-4&0&-1&4&-4&0 \\
  0&0&0&0&1&2&0&0&-1&2&0 \\
  0&0&0&0&0&0,5&1&0&0&-0,5&1 \\
  \hline
  \end{tabular}

  Table 6

   \begin{tabular}{|c|c|c|c|c|c|c|c|c|c|c|}
  \hline
  \hline
  1&2 & 3 & 4 &5&6&7&8& 9&10&11 \\
  \hline
  100&0&-100&0&-8,75&-100&0&-5&3,75&0&-575 \\
   1&1&1&0&0,25&1&0&0&-0,25&0&2 \\
  0&0&0&1&2&0&0&-1&2&0&0 \\
  0&0&0&0&0,5&1&0&0&-0,5&1&0 \\
  0&0&0&0&0,25&1&1&0&-0,25&0&1 \\
  \hline
  \end{tabular}

  Table 7

   \begin{tabular}{|c|c|c|c|c|c|c|c|c|c|c|}
  \hline
  \hline
  1&2 & 3 & 4 &5&6&7&8& 9&10&11 \\
  \hline
  100&0&-100&-1,875&-12,5&-100&0&-3,125&0&0&-575 \\
   1&1&1&0,125&0,5&1&0&-0,125&0&0&2 \\
  0&0&0&0,5&1&0&0&-0,5&1&0&0 \\
  0&0&0&0,25&1&1&0&-0,25&0&1&0 \\
  0&0&0&0,125&0,5&1&1&-0.125&0&0&1 \\
  \hline
  \end{tabular}

  Table 8

Table 5 shows the solution of the problem for $-z\geq 500, x=(5, 5, 65/2, 375)^T$ . Next, we continue according to the second column, to which the inequality $-z \geq 600$ initially corresponded. Based on the first row of the transformed table 6, the pivoting column is 10 and in the next step, column 9. Table 8 shows the solution to the system in case of constraint $-z \geq 600, x= (3,125, 0, 0,575)^T$. In Table 9, it is in case of constraint $-z\geq 700$, the objective function of the dual problem (3) is unbounded (column 8), the system is contradictory.

    \begin{tabular}{|c|c|c|c|c|c|c|c|c|c|c|}
  \hline
  \hline
  1&2 & 3 & 4 &5&6&7&8& 9&10&11 \\
  \hline
  0&-100&-200&-14,375&-62,5&-200&0&9,375&0&0&-775 \\
   1&1&1&0,125&0,5&1&0&-0,125&0&0&2 \\
  0&0&0&0,5&1&0&0&-0,5&1&0&0 \\
  0&0&0&0,25&1&1&0&-0,25&0&1&0 \\
  0&0&0&0,125&0,5&1&1&-0.125&0&0&1 \\
  \hline
  \end{tabular}

  Table 9

In fact, you don't need three columns for the objective function, just one. If we have already found a feasible solution and we want an even smaller objective function value, we put a positive number instead of zero in the objective function column of the w-row. Otherwise, the zero in the w-row must be replaced by a negative number, narrowing the objective function less. Using various systems of inequalities, we determine a solution in which the objective function has a value that satisfies us, has an appropriate number of non-zero components of the vector $x$, etc. Based on the calculations, we can also specify the setting of the problem, change constraints, add variables, etc.

    If the objective function of the problem can be unbounded, then a column with coefficients of the objective function can be added to the table, where in the w-row is a sufficiently large number (700 in our example).

    Let us consider a linear programming problem

    \begin{equation}\
\ z= (c,x) \rightarrow min,
\ {Ax \geq b.}
\ x\geq 0,
\end{equation}
 and its dual
    \begin{equation}\
\ w = (y,b) \rightarrow max,
\ {yA \leq c.}
\ y\geq 0,
\end{equation}
where $y$ is a $m-$ vector of dual variables.

    If the dimension of the problem is not large, we replace them with an equivalent system of inequalities

   \begin{equation}\
\ -(c,x)+(y,b) \geq 0,
\end{equation}
   \begin{equation}\
\ {Ax\,\,\,\,\,\, \,\,\,\,\,\,\,\,\,\,\,\,\,\,\,\,\,\,\geq b.}
\end{equation}
   \begin{equation}\
\\,\,\,\,\,\,\,\,\,\,\,\, \,\,\,\,\,\,\,\,\,\,\,\,-yA \geq -c,
\end{equation}
   \begin{equation}\
\ x\geq 0,y\geq 0.
\end{equation}

The system of inequalities (6)-(9) composed for Example 3.2 was solved in three steps, $ x=(0; 2; 0; 2)^{T} , y = (0, 0,1).$

\vskip 1.5cm
\section{\bf Conclusion}

Solved examples of linear inequalities demonstrate the effectiveness of the presented method. The number of steps required almost never exceeded the number of inequalities.

   \vskip 1.0cm
    \textbf{Citations}
     \vskip 1.0cm

    [1] L. Khachiyan, \textit{Fourier-Motzkin elimination method}, Encyklopedia Of Optimization, (2001) vol.2,155-159.

    [2]  E. Ubi, \textit{Mathematical programming via the least-squares method}, Cent. Eur. J. Math., vol. 8, (2010), 795-806.

    [3]  E.Ubi,  \textit{Linear inequalities via least squares}, Proc. of the Estonian Academy of Sciences, (2013), 62, 4, 238-248.

    [4] T.Hu, \textit{Integer programming and network flows}, Addison-Wesley Publishing Company, Massachusetts, (1970),519 pp.

    [5] S. Chernikov, \textit{ Lineare Ungleichungen} , Deutscher Verlag der Wissenschaften, Berlin, (1971),486 pp.

    [6]  S.Gass and S.Vinjamuri:  \textit{Cycling in Linear Programming Problems}, Comp.Oper. Res.,31(2004), 303 - 311.

    [7] C.H.Papadimitriou and K.Steiglitz: \textit{Combinatorial Optimization: Algorithms and Complexity}, Dover Publications, New York (1998), 446 pp.

\end{document}